%% file: rh.tex
\documentclass[12pt]{article}
\usepackage{graphics}

\date{}
\def\qed{{\hfill\rule{1.2ex}{1.2ex}}}
\newtheorem{Th}{Theorem}
\newtheorem{Cor}{Corollary}
\newenvironment{Proof}{\begingroup\noindent{\bf Proof:}}{\qed\medskip\endgroup}

\newtheorem{Lemma}{Lemma}

\def\dim{\mathop{\rm dim}}

\def\C{{\bf C}}
\def\R{{\bf R}}
\def\Q{{\bf Q}}
\def\N{{\bf N}}
\def\Z{{\bf Z}}

\author{Norbert A'Campo}
\title{Monodromy of real isolated singularities} 

\begin{document}

\maketitle

\begin{abstract}{Complex conjugation on complex space permutes the level sets 
of a real polynomial function and induces involutions 
on level
sets corresponding to
real values. For isolated complex hypersurface singularities with 
real defining equation we show the existence of a 
monodromy vector field such that 
complex conjugation intertwines 
the local monodromy  diffeomorphism  
with its inverse. In particular, 
it follows that the geometric
monodromy is the composition of the involution induced by
complex conjugation and another involution. This topological property 
holds for all isolated complex plane curve singularities. Using real 
morsifications, we compute the action of
complex conjugation and of the other involution on the 
Milnor fiber of real plane curve singularities.\\[1mm]
}
\end{abstract}

\noindent
{\bf Introduction.}
Let $f:\C^{n+1} \to \C$ be a map defined by a 
polynomial. We assume that $f(0)=0$ and that $0 \in
\C^{n+1}$ is an isolated critical point of $f$. For $p \in \C^{n+1}$ let 
$||p||$ denote the square root of $|z_0(p)|^2+|z_1(p)|^2+ \cdots +|z_n(p)|^2$.
Let $B_{\epsilon}:=\{p\in \C^{n+1} \mid ||p|| \leq \epsilon \}, 0 <\epsilon,$
be a Milnor ball for the singularity of $f$ 
and let 
${\rm Tube}_{\epsilon,\delta}:=\{ p \in B_{\epsilon} \mid |f(p)| \leq \delta \}, 0 < \delta << \epsilon,$
be a regular tubular neighborhood of $\{p \in B_{\epsilon} \mid f(p)=0 \}$ in
$B_{\epsilon}$. A monodromy vector field $X$ for the singularity is a smooth 
vector field
$p \in {\rm Tube}_{\epsilon,\delta} \mapsto X_p \in \C^{n+1}$ such 
that we have the following properties for 
$p \in {\rm Tube}_{\epsilon,\delta}$ 
(remember $i=\sqrt{-1}$): 

- $(df)_p(X_p)=2\pi i f(p)$, 

- $X_p$ is tangent to $\partial{B_{\epsilon}}$ if 
$p \in \partial{B_{\epsilon}}$,

- trajectories of $X$  starting at $p \in \partial{B_{\epsilon}}$
are periodic with period $1$ and are the boundary of a smooth disk in 
$\partial{B_{\epsilon}}$,
that is transversal to the function $f$.

Using partition of unity, one can construct monodromy vector fields. The
flow at time $1$ of a monodromy vector field $X$ defines a
monodromy diffeomorphism 
$T=T_X\colon F \to F$, where the manifold with boundary
$(F,\partial(F)):=\{ p \in B_{\epsilon} \mid f(p)=\delta \}$ is 
the Milnor fiber of the singularity. The 
relative isotopy class 
of the diffeomorphism $T$ is independent from the chosen 
monodromy vector field and is the geometric monodromy of the singularity. 
The geometric monodromy is a topological invariant of the singularity.

From now on we will assume in addition, that the polynomial $f$ is real 
meaning that its coefficients are real numbers. Let $c:\C^{n+1} \to \C^{n+1}$
denote the involution on complex space 
given by the complex conjugation of coordinate values.
Hence with the above notations, we 
have $c({\rm Tube}_{\epsilon,\delta})={\rm Tube}_{\epsilon,\delta}$ and 
$c(F)=F$. We denote by $c_F:F \to F$ the restriction of the involution 
$c$ to $F$.

Let $X\colon {\rm Tube}_{\epsilon,\delta} \to \C^{n+1}$ be a 
monodromy vector field for the isolated singularity of $f$. We may 
assume that we have 
constructed the vector field $X$ with more care near the boundary of 
the Milnor ball in order to achieve  
that for some 
$\epsilon' < \epsilon$
we have the symmetry
$c(X_p)=-c(X_{c(p)}), p \in {\rm Tube}_{\epsilon,\delta}, ||p|| > \epsilon'$.

\noindent
Since $f$ is real, we have
$$
c((df)_p(X_p))=(df)_{c(p)}(c(X_p))=-2\pi i f(c(p))
$$
hence, we see (by substituting $q$ for $c(p)$ and accordingly $c(q)$ for $p$) 
that the vector field $X^c$ 
defined by:
$$
q \in {\rm Tube}_{\epsilon,\delta} \mapsto X^c_q:=-c(X_{c(q)}) \in \C^{n+1},
$$
is  a monodromy vector field too. Let 
$Y:{\rm Tube}_{\epsilon,\delta} \to \C^{n+1}$ be the  vector field
$Y:={X+X^c \over 2}$, which due to the extra care is also 
a monodromy 
vector field. We have $Y^c=Y$. 
The following is an important symmetry of the
geometric monodromy:

\begin{Lemma} Let $T_Y$ be a  monodromy diffeomorphism, which has 
been computed with a monodromy vector field $Y$ satisfiyng  
$Y^c=Y$. We have the symmetry
$$
c_F \circ T_Y \circ c_F =T_Y^{-1}.
$$ 
The geometric
monodromy $T$ satisfies (up to relative isotopy) the symmetry
$$
c_F \circ T \circ c_F =T^{-1}.
$$
\end{Lemma}

\begin{Proof} The restriction 
of complex conjugation 
$c_{{\rm Tube}_{\epsilon,\delta}}\colon {\rm Tube}_{\epsilon,\delta} 
\to {\rm Tube}_{\epsilon,\delta}$ maps the monodromy vector field
$Y$ to $-Y$ and $F$ to $F$. Hence, since $c_{{\rm Tube}_{\epsilon,\delta}}$ 
reverses the orientations of the 
trajectories, we have 
$$
T_Y^{-1}=(c_F)^{-1} \circ T_Y \circ c_F=c_F \circ T_Y \circ c_F.
$$
Since the geometric monodromy $T$ is in  the  relative mapping class 
group of the Milnor fiber represented by $T_Y$, for $T$ in the
relative
mapping class group we have the symmetry $c_F \circ T \circ c_F =T^{-1}.$
\end{Proof}

Symmetries of monodromies as in the lemma 
can occur in the more general context of  so-called 
{\it strongly invertible knots}, see for instance
[H-To],[T],[K]. 

The symmetry property $c_F \circ T \circ c_F =T^{-1}$ expresses that
the geometric monodromy $T$ of a complex hypersurface with real 
defining equation is conjugate in the mapping class group by an 
element of order $2$ to its 
inverse $T^{-1}$. This statement  
does not refer to  any complex conjugation, so it can be stated for
any complex hypersurface singularity with complex defining equation. We say 
that the singularity  is {\it strongly invertible} if its geometric monodromy
diffeomorphism
$T$ is conjugate by an
element of order $2$ 
in the relative mapping class group of the Milnor fiber 
to its
inverse $T^{-1}$. 
The property of strong invertibility is a topological property for hypersurface
singularities.

We can rewrite the symmetry property as follows:
$T_Y \circ c_F \circ T_Y \circ c_F= {\rm Id}_F$. We see that 
$T_Y \circ c_F:F \to F$ is an involution of $F$. It follows the 

\begin{Cor} The geometric monodromy $T$ of an 
isolated  complex hypersurface singularity, which is 
defined by a real equation, is the composition of two involutions of the 
fiber $F_\delta, \delta \in \R,$ namely: 
$T= (T \circ c_F) \circ c_F$, where $c_F$ is the restriction of the complex 
conjugation.
\end{Cor}

For $k \in \Z$ we also have the
relation
$T_Y^k \circ c_F \circ T_Y^k \circ c_F= {\rm Id}_F$,
which shows that 
$T_Y^k \circ c_F:F \to F, k \in \Z,$ is a sequence of involutions of $F$.

The above observations can be applied to plane curve singularities 
in general, since
it follows from the Theory of Puiseux Pairs that
every plane curve singularity is topologically equivalent to a singularity
given by a real equation and hence plane curve singularities are 
strongly invertible.

For complex hyper\-surface singu\-larities of higher dimension the situation
seems to be  op\-posite.
In $\C^{n+1},\,n>1,$ there exist iso\-lated hyper\-surface singula\-rities  which are
not topologically equivalent to a 
singularity with a real defining equation [T].
We expect that in general the geometric monodromy $T$ 
of a complex hypersurface singularity fails to be strongly invertible.  

For  complex hypersurface singularities
the eigenvalues of the monodromy $T_*$ acting on
the homology $E:=\bigoplus_k H_k(F,\Q)$  are
roots of unity by the monodromy theorem and  the characteristic polynomial
of $T_*$ is a product of cyclotomic polynomials. It follows, that the complex homological monodromy ist
strongly invertible. Mathias Schulze has proved strong invertibility of the homological monodromy
with real coefficients [S] and the following is a strengthening of his
result.

\begin{Th}
The rational homological monodromy of a complex hypersurface singularity 
is strongly invertible. 
\end{Th}

\begin{Proof}
Let
$E=\bigoplus_i E_i$ a be finest possible direct sum 
decomposition in $T_*$-invariant 
$\Q$-subspaces of $E$. The characteristic polynomial of the restriction
$A_i:E_i \to E_i$ 
of $T_*$ to a summand $E_i$ is a power $\psi(t)=\phi(t)^L$ of a cyclotomic 
polynomial $\phi(t)$. We have $\psi(A_i)=0$ by the Hamilton-Cayley 
theorem and we have $\psi(A_i^{-1})=0$ since a power of cyclotomic polynomial
satisfies $\psi(1/t)=\pm t^{{\rm degree}(\psi)} \psi(t)$. 

Since the decomposition has no refinement, we may
choose a  vector $e_1 \in E_i$,  that is 
cyclic for $A_i$ and for $A_i^{-1}$. 
The systems
$$
e_1,e_2:=A_i(e_1), \cdots ,e_{\dim E_i}:=A^{\dim E_i -1}(e_1),
$$
$$
f_1:=e_1,f_2:=A_i^{-1}(e_1), \cdots ,f_{\dim E_i}:=A^{-\dim E_i +1}(e_1)
$$ 
are bases
for the space $E_i$. Let $b_i:E_i \to E_i$ be the linear map defined by
$b_i(e_j)=f_j, 1 \leq j \leq \dim E_i$. 

We have
$A_i^{-1}b_i(e_j)=b_iA_i(e_j), 1 \leq j < \dim E_i$. The polynomial 
$\psi_1(t):=-\psi(t)+t^{\dim E_i}$ satisfies 
$\psi_1(A_i^{\pm 1})=A_i^{\pm \dim E_i}$. Since the degree of 
the polynomial
$\psi_1(t)$  is less than $\dim E_i$, we have
$$
A_i^{-1}b_i(e_{\dim E_i})=A^{-\dim E_i}b_i(e_1)=\psi_1(A_i^{-1})b_i(e_1)=
$$
$$
b_i\psi_1(A_i)(e_1)=b_iA_i^{\dim E_i}(e_1)=b_iA_i(e_{\dim E_i}).
$$
We conclude that $A_i^{-1}b_i=b_iA_i$ and $A_i^{-1}=b_iA_ib_i^{-1}$ hold.

The polynomial
$\psi_0(t):={ {\psi(t)-\psi(0)} \over {-t\psi(0)} }$ satisfies
$A_i^{-1}=\psi_0(A_i)$ and $A_i=\psi_0(A_i^{-1})$. 
We deduce $b_iA_i^{-1}=b_i\psi_0(A_i)=\psi_0(A_i^{-1})b_i=A_ib_i$ and 
conclude $A_i^{-1}=b_i^{-1}A_ib_i$. We observe at this point that both 
the conjugates of $A_i$ by $b_i$ and by
$b_i^{-1}$ are equal to the inverse $A_i^{-1}$.

For $0 \leq j < \dim E_i$ we have (remember $e_1=f_1=b_i(e_1)=b_i(f_1)$)
$$
b_i(f_{j+1})=b_iA_i^{-j}(f_1)=A_i{^j}b_i(f_1)=A_i{^j}(e_1)=e_{j+1}.
$$
Hence
$b_i$ is of order two, which shows that $A_i$ is strongly invertible over
$\Q$. The sum $b:=\bigoplus b_i$ is a rational strong inversion  for $T_*$. 
\end{Proof}

\noindent
{\bf Question:} Is the integral homological monodromy of a complex hypersurface singularity
strongly invertible? We would like to know the answer for the isolated surface singularity
$f(x,y,z)=L_1 L_2^2 L_3^3 L_4^4 L_5^5 L_6^6 L_7^7 L_8^8 L_9^9 +x^{46}+y^{46}+z^{46}$
where $L_j, 1 \leq j \leq 9$ are linear forms on $\C^3$, such that the $9$ lines $\{L_j=0\}$ in the complex 
projective plane
span the $9_3$ configuration of flex tangents to a nonsingular cubic. No real equation for this singularity
can exist, since the configuration $9_3$ cannot be realized in the real projective plane, see [H-C,T].

In Section {\bf 1} we will study in detail the effect of complex
conjugation on the topology of plane curve singularities. We 
show that the fiber of the link of a
connected divide carries naturally a cellular decomposition with tri-valent
$1$-skeleton.

In Section {\bf 2} we study the involutions that appear in the decomposition
of the geometric monodromy of plane curve singularities. We show that they lift to $\pi$-rotations
about an axis in the universal cyclic covering of the complement of the link.

I like to thank Makoto Sakuma for explaining  the work
of Jeffrey L. Tollefson [To] to me.

I like to thank the referee for having attracted my attention to work of Sabir Gusein-Zade [G-Z2], in which he shows
among other results, that the integral homological monodromy of an isolated complex hypersurface singularity
with real defining equation $f$ is the composition of two involutions, both conjugated to the action of complex
conjugation on the homology of a real regular fiber.

\noindent
{\bf 1. Real plane curve singularities.} 
Let $f:\C^2 \to \C$ define an isolated plane 
curve singularity  at $0 \in \C^2$ given by a real convergent power series  
$f\in \R\{x,y\}$. Let $f=f_1f_2 \cdots f_r$ be the factorization in
local branches.  A factor $f_i\in \C\{x,y\}$ is a 
real convergent power series, i.e. $f_i(\R^2) \subset \R$ or the conjugate
series $c(f_i)$ is a factor too. 

The topology of a plane curve singularity $f$ is completely encoded in 
a divide for $f$ [AC3,B-K], see also [AC1-2,G-Z1]. 

We state the results  more generally for links of 
connected divides, since the involution given by complex conjugation 
on real isolated plane curve singularities 
corresponds to a natural involution of
links of divides as explained below, see [AC3-4].

A divide is the image of a generic relative immersion $P$
of a compact one dimensional manifold in the unit disk $D$ in $\R^2$. The 
complement in $S^3$ of the link
$L(P)$
of $P$ is naturally fibered over $S^1$, if the divide $P$ 
is connected. 
Let $\delta_P$ be the number of double 
points of the divide $P$. We denote the  fiber surface 
above $1$ by $F_P$.
The local topology of a plane curve singularity is obtained from
a divide for the singularity.  More precisely, 
the link $L_f \subset \partial{B_{\epsilon}}$ of an isolated 
plane curve singularity
$\{f=0\}, f \in \R\{x,y\},$
is equivalent to the link  $L_P$ of a divide $P$, see [AC3].  

We recall, 
that for a plane curve singularity $f$ one can obtain
a divide $P$ by performing 
a small real deformation $f_s, 0 \leq s \leq 1,$ of the 
singularity, where for $0 < s \leq 1$ the restriction of $f_s$ to 
the euclidean disk $D_{\epsilon}:=B_{\epsilon} \cap \R^2$ 
of radius $\epsilon$ in $\R^2$ is a Morse function with 
$\mu(f)$ critical points and such
that the $0$-level is connected and contains all the saddle points. 
The divide $P$ for $f$ is the curve $\{ p \in D_{\epsilon} \mid f_1(p)=0 \}$, 
which we rescale by the factor ${1 \over \epsilon}$ from 
$D_{\epsilon}$ into the 
unit disk $D$.

\begin{center}
\input{fig-1.pstex_t}
\newline
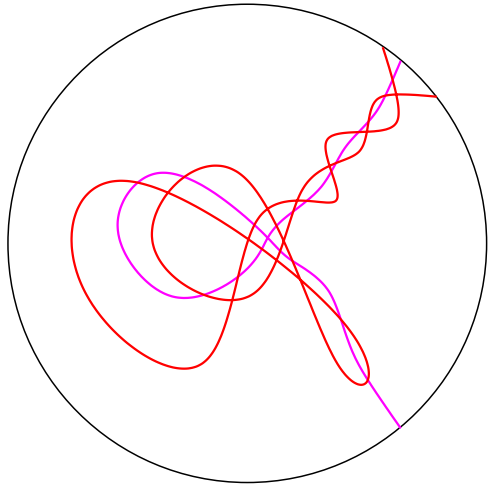
{Fig. 1. Divide for the singularity $(x^3-y^2)((x^3-y^2)^2-4x^8y)$.}
\end{center}

The homology $H_1(F_P,\Z)$ can be described combinatorially in terms of the 
divide $P$ as a direct sum 
$H_1(F_P,\Z)=E_- \bigoplus E_0 \bigoplus  E_+$, where
$E_-,E_0$ and $E_+$ are the subspaces in 
$H_1(F,\Z)$, which are freely generated  
as follows: 
$E_+:=[\delta_1, \cdots ,\delta_{\mu_+}]$, 
$E_0:=[\delta_{\mu_++1}, \cdots ,\delta_{\mu_+ + \mu_{0} }]$ 
and 
$E_-:=[\delta_{\mu_++\mu_{0}+1}, \cdots , \delta_{\mu_-+\mu_{0}+\mu_+}]$, 
where $(\delta_i)_{ 1 \leq i \leq \mu}$ is 
the oriented system of vanishing cycles of the 
divide $P$ with positive upper triangular Seifert form 
$S:H_1(F,\Z) \to H^1(F,\Z)={\rm Hom}(H_1(F,\Z),\Z)$, 
see [AC4] Theorem $4$. We define $N:=S - {\rm Id}$, which is 
upper-triangular nilpotent matrix. The monodromy
$T_*:H_1(F,\Z) \to H_1(F,\Z)$ is given by:
$T_*=(S^t)^{-1} \circ S.$ 

We recall briefly the construction of the 
fiber 
surface  $F_P$ of a divide $P$, see [AC4]. 
Let $P$ be a connected divide in the unit euclidian 
disk $D$ consisting of a generic system of 
immersed copies of $S^1$ and $[0,1]$. We assume, that $P$ is 
irreducible, i.e. every differentiable 
chord of
$D$ transversal to $P$, that meets $P$ in at 
most $1$ point, is transversally 
isotopic to an arc in the boundary of $D$. Let 
$f:D \to \R$ be a Morse function adapted to $P$. A region of $P$ 
is a connected component of $D \setminus P$, that does not meet the 
boundary of $D$. If $P$ 
has regions, we assume that the sign of $f$ is chosen such that in 
at least one region the 
function $f$ is positive. The link $L_P$ of $P$ is the subset
in the $3$-sphere sitting in the  tangent space 
$S^3:=\{ (x,u) \in T\R^2 \mid ||x||^2+||u||^2=1 \}$ and is given by
$L_P:=\{(x,u) \in S^3 \mid u \in TP \}$. 
The number $r=r(P)$ 
of components of the link
$L_P$ is twice the number of immersed circles plus 
the number immersed intervals in $P$.

\begin{center}
\input{fig-2.pstex_t}
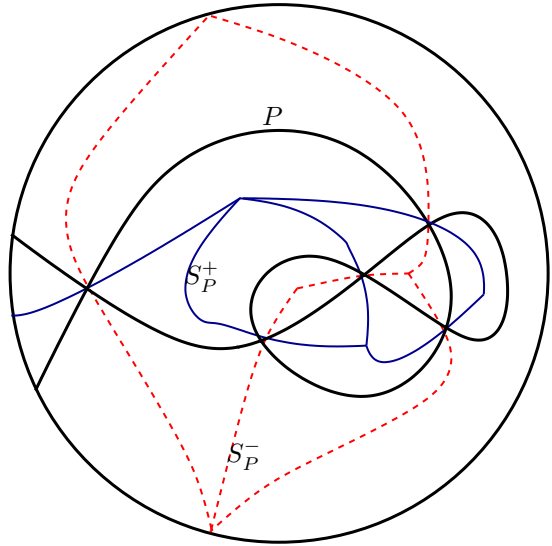
\center{Fig. 2. Divide and graph of divide.}
\end{center}

Let $S_P$ be a system
of gradient lines of $f$, which connect saddle points of $f$ with  
local or relative maxima or minima of $f$. In fact, near the boundary of 
$D$ the lines of the system $S_P$ are only gradient-like, but end in a 
relative critical point of $f$. We call the  
set $S_P$ the graph of the divide $P$. By a slalom 
construction, see [AC5], one can reconstruct
from each of the sets $S_P^{\pm}:=\{x \in S_P \mid \pm f(x) \geq 0 \}$ 
the divide $P$. We assume that at the 
maxima and minima of $f$, that different
arcs of $S_P$ have different unoriented tangent directions. Let $\Sigma_P$
be the following subset of the fiber $F_P$ over $1$ of the natural fibration
over $S^1$ of the complement of the link of $P$: the subset $\Sigma_P$ is
the closure in the tangent space $TD$ of $D$ of
$$
\Sigma'_P:=\{(x,u) \in TD \mid x \in S_P, f(x)> 0,
(df)_x(u)=0, ||x||^2+||u||^2=1 \}.
$$

\noindent
We recall that the union of fiber $F_P$ with the link $L_P$
is the closure in $TD$ of
$$
F'_P:=\{(x,u) \in TD \mid x \in D, f(x) > 0,
(df)_x(u)=0, ||x||^2+||u||^2=1 \},
$$
\noindent
which is a surface of genus $g(P)$ equal to the number of double points in $P$.

\begin{center}
\input{fig-3.pstex_t}
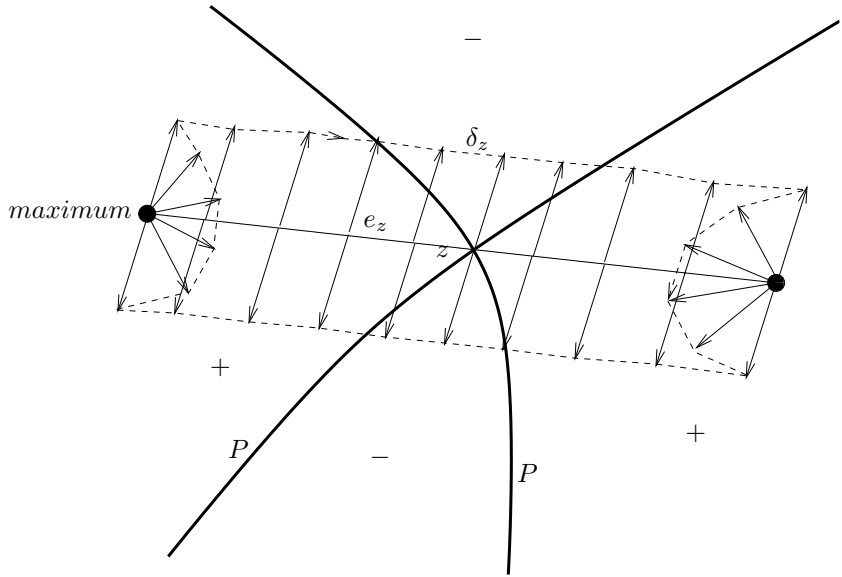
\center{Fig. 3. Vanishing cycle $\delta_z$ for a saddle point.}
\end{center}

The system of cycles $\delta_i, 1 \leq i \leq \mu,$ can be drawn on 
$\Sigma_P$ as follows. The vanishing cycle 
$\delta_i, \mu_-+\mu_0 < i \leq \mu,$
which corresponds to a local maximum $M$ of $f$ is the circle
$\{ (M,u) \mid u \in T_MD, ||M||^2+||u||^2=1\} $ oriented counter clock-wise.
The vanishing cycle $\delta_i, \mu_- < i \leq \mu_-+\mu_0,$ which 
corresponds to a saddle point $z$ of $f$ is a curve in the set
$E_z:=\{(x,u) \in TD \mid x \in e_z,
u\in {\rm Kernel}(df)_{e_z}, ||x||^2+||u||^2=1 \}$, as drawn in Fig-3. 
The curve $\delta_i$ is a piecewise smooth embedded  
copy of $S^1$ in $F_P$ with image in $E_z$
and with non-constant projection to the disk $D$. The orientation 
is chosen such that at both ends the orientation agrees with the oriented
vanishing cycle of the maximum. Moreover, the inward tangent vectors
$(x,u)$ to $e_z$ at end points $x \in e_z$ do  belong to $\delta_i$. 
The vanishing cycle $\delta_i, 1 \leq  i \leq \mu_-,$ which 
corresponds to a local minimum of $f$ projects to an oriented circuit
$e_1,e_2, \cdots ,e_k$ of edges of the graph $S_P$. The circuit 
surrounds the region in counter clock-sense 
to which 
the minimum corresponds. The circuit is a polygon and bounds a cell in $D$.
The curve $\delta_i$ is a subset of
$\{(x,u) \in TD \mid x \in \cup_{1 \leq j \leq k}  e_j, u\in {\rm Kernel}(df)_{x}, ||x||^2+||u||^2=1 \}$
and is the image of a piecewise smooth  embedding of 
$S^1$ as drawn in Fig-4. The vectors
$(x,u)$, which belong to $\delta_i$, point out of the cell of the circuit. 
Smooth representatives for the system of vanishing cycles can be obtained 
using tears as in [AC6].

\begin{center}
\input{fig-4.pstex_t}
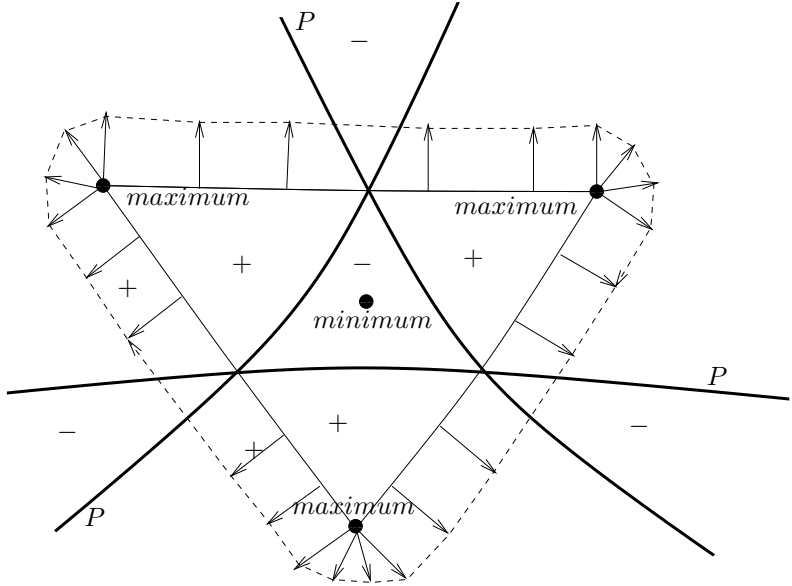
\center{Fig. 4. Vanishing cycle $\delta_i$ for a minimum.}
\end{center}

The involution
$c_{T\R^2}:(x,u)\in T\R^2 \mapsto (x,-u)\in T\R^2$
induces an involution
$c_{S^3}$ on $S^3\subset T\R^2$ that preserves the link $L_P$
of any divide $P$ and that
induces involutions on the fiber surfaces above $\pm 1$.
If the divide $P$ is a divide for a real plane curve singularity $f$, the
involutions induced by $c_{T\R^2}$ or by the complex conjugation $c$ on
the triples $(S^3,L_P,F_P)$ and on $(\partial{B_{\epsilon}},L_f,F)$ 
correspond to
each other by the homeomorphism of pairs of the main theorem of [AC3].

\begin{Th} Let $P \subset D$ be a connected divide consisting of the 
image of a relative generic immersion of a compact one dimensional manifold   
in the unit disk 
$D \subset \R^2$. Let $E_-,E_0$ and $E_+$ be the summands in 
$H_1(F,\Z)$ as above. 
The involution $c_{F*}:H_1(F,\Z) \to H_1(F,\Z)$ fixes pointwise
the summand $E_+$, in particular for $\delta_i \in E_+$ we have
$$
c_{F*}(\delta_i)=\delta_i.
$$
For $\delta_i \in E_0$ we have
$$
c_{F*}(\delta_i)= -\delta_i+
\sum_{1\leq j \leq \mu_+} \langle N(\delta_j),\delta_i \rangle.
$$
For $\delta_i \in E_-$ we have
$$
c_{F*}(\delta_i)= \delta_i+
\sum_{1\leq j \leq \mu_+} \langle N(\delta_j),\delta_i \rangle-
\sum_{\mu_+ < j \leq \mu_++\mu_0} 
\langle N(\delta_j), \delta_i \rangle \delta_j.
$$
The trace  
of the involution $c_{F*}:H_1(F,\Z) \to H_1(F,\Z)$
is given by: 
$$
{\rm Trace}(c_{F*})=\mu_- - \mu_0 + \mu_+.
$$
\end{Th}  

\begin{Proof} A vanishing cycle $\delta_i \in E_+$ corresponds 
to a maximum $M \in D$ of $f$,
hence as set we have $\delta_i=\{(M,u) \in TD \mid ||M||^2+||u||^2=1\}$.
The
involution $(x,u) \mapsto (x,-u)$ induces on $\delta_i$ the antipodal map,
which is orientation preserving. It follows $c_{F*}(\delta_i)=\delta_i$ in
homology.  
A vanishing cycle $\delta_i \in E_0$ corresponds to a saddle point 
$z \in D$ of $f$. Working with the tear model  of [AC6], see also Section 2, 
we see that the involution reverses the orientation and that at
the endpoints, which are maxima of $f$, we have outward instead of 
inward vectors. Hence, 
$$
c_{F*}(\delta_i)=-\delta_i+\sum m_{i,j} \delta_j,
$$ 
where in the sum
$j$ runs through the maxima of $f$ in the interior of $D$ which are 
connected by gradient 
lines of $S_P$ to the saddle point $i$. The coefficient $m_{i,j}$ equals $1$ or $2$ 
depending on whether the connection by gradient lines is simple or double.
Finally one gets 
$c_{F*}(\delta_i)= -\delta_i+
\sum_{1\leq j \leq \mu_+} \langle N(\delta_j),\delta_i \rangle$. 
For a vanishing cycle $\delta_i \in E_-$
one can work with the model of [AC6]  and get 
$c_{F*}(\delta_i)= \delta_i+
\sum_{1\leq j \leq \mu_+} \langle N(\delta_j),\delta_i \rangle-
\sum_{\mu_+ < j \leq \mu_+ \mu_0} \langle N(\delta_j),\delta_i \rangle$. 
Since $N$ is strictly upper-triangular, the matrix of $c_{F*}$ is 
upper-triangular with $\pm 1$ on the diagonal. We get 
${\rm Trace}(c_{F*})=\mu_- - \mu_0 + \mu_+$.
\end{Proof}

We can present the Seifert form and homological 
monodromy with block-matrices $S,T$ as in [G-Z2]. On $H_1(F,\Z),E_-,E_0,E_+$ and 
$H^1(F,\Z)$ we work with the basis or dual basis given by the system 
$\delta_i, 1 \leq i \leq \mu$. 
The matrix in block form of the Seifert form is
$$
S=\left(\matrix{{\rm Id}_{\mu_+}& A & G\cr O& {\rm Id}_{\mu_0}& B\cr
O & O & {\rm Id}_{\mu_-}\cr}\right)  
$$
where the block $G$ equals the block matrix product $1/2 (A \circ B)$. 
The matrix coefficients of $A \circ B$ and $G$ have interpretations 
in terms of the divide $P$ or in terms of the Morse function 
$f_P$ on the disk $D$. 
The matrix coefficient 
$(A \circ B)_{(i,j)}, 1 \leq i \leq \mu_+, \mu_+ + \mu_0 < j \leq \mu,$ 
counts the number of  sector adjacencies, that has the 
$+$-region $i$ of $P$  with the $-$-region $j$, while the 
coefficient $G_{(i,j)}$ counts the number of 
common boundary 
segments of the regions $i$ and $j$. The coefficient $(A \circ B)_{(i,j)}$
also counts the number of 
saddle connections via gradient lines of $f_P$ in between
the minimum $j$ and maximum $i$, while the coefficient $G_{(i,j)}$
counts the number of components of regular gradient line connections from
the minimum $j$ to the maximum $i$.
This explains the above factor $1/2$, since a segment of 
$P$ is twice incident 
with a saddle point of $f_P$. 

The matrix of the action of complex conjugation on $H_1(F,\Z)$ is
$$
C=\left(\matrix{{\rm Id}_{\mu_+}& A & G\cr 
O& {\rm -Id}_{\mu_0}& -B\cr
O & O & {\rm Id}_{\mu_-}\cr}\right)  
$$
and is obtained from the matrix $S$ by multiplying the middle row of blocks by $-1$.
It is interesting to compute the matrix of the involution $T\circ c_F$ 
on $H_1(F,\Z)$
$$
TC=\left(\matrix{{\rm Id}_{\mu_+}& O & O\cr 
-{^tA}& {\rm -Id}_{\mu_0}& O \cr
{^tG} & {^tB} & {\rm Id}_{\mu_-}\cr}\right)  
$$
The matrix $T\circ C$ is the transgradient of the matrix of $C$. 

It turns out that 
the combinatorial property $G=1/2 (A \circ B)$ for divides 
is equivalent to 
$C\circ C={\rm Id}_{\mu}$ or $T\circ C\circ T\circ C={\rm Id}_{\mu}$.

For an isolated  plane curve singularity at $0 \in \C^2$, which is 
given by a real 
equation $\{f=0\}, f \in \R\{x,y\}, f(0)=0,$ we will 
denote by $\delta_{\R}(f)$ the number of 
double points of a divide $P$ 
for the singularity.  Hence $\delta_{\R}(f)$ is 
the maximal number of local
real 
saddle points in some level, that can occur for  
a small real deformation of $f$.
Observe, that one has 
$\delta_{\R}(f) \leq \delta(f)$, where $\delta(f)$ 
is the maximal number of local
critical points in some level, that can occur for  
a small deformation of $f$. We recall the formula
$\mu(f)=2 \delta(f)-r+1$ of Milnor [M]. As example for 
$f=x^4+Kx^2y^2+y^4, -2 \not= K \not= 2,$ one has
$\delta_{\R}(f)=4, \delta(f)=6, r=4$ and 
very surprisingly, Callagan shows 
that
for $-2 < K < 2$ there exists a small real deformation with
$5$ local minima in the same level [C].

\begin{Th} For an isolated  plane curve singularity at 
$0 \in \C^2$, which is  given by a real 
equation $\{f=0\}, f \in \R\{x,y\}, f(0)=0,$ we have
$$
\mu(f)=2 \delta_{\R}(f)+{\rm Trace}(c_{F*})
$$
\end{Th}

\begin{Proof} 
We have $\delta_{\R}={\rm dim}E_0=\mu_0.$
Hence 
$\mu(f)=\mu_-+\mu_0+\mu_+=
{\rm Trace}(c_{F*})+2\mu_0=2 \delta_{\R}(f)+{\rm Trace}(c_{F*}).$
\end{Proof}

Curves $\delta_i$ that correspond to maxima of $f_P$ are invariant by the involution $c$. A curve $\delta_i$ 
that corresponds to a saddle point or minimum of $f_P$ is in general not invariant by $c$. The union
$\Sigma_P$ of the curves $\delta_i$ is invariant by $c$. We get

\begin{Th} Let $P$ be a connected divide. The pair $(F_P,\Sigma_P)$ 
defines  a tri-valent cellular decomposition with $r(P)$ punctured 
cells  
of the fiber surface with boundary $F_P$. The involution 
$c:(x,u) \to (x,-u)$ 
acts on the map $(F_P,\Sigma_P)$.
\end{Th}

\begin{Proof}
The union $\Sigma_P$ of the vanishing cycles $\delta_i$
is a tri-valent graph $\Sigma_P$ that is invariant by $c$. The inclusion $\Sigma_P \subset F_P$ induces an isomorphism 
$H_1(\Sigma_P,\Z) \to H_1(F_P,\Z)$, hence, $\Sigma_P$ is a spine for the surface $F_P$.
\end{Proof}

In particular if $r(P)=1$ the triple $(F_P,\Sigma_P,c)$ is a so-called
maximal unicell map of genus $g(P)$ with orientation reversing 
involution $c$ with $r(P)=1$ fixed points on the graph of the map. 
The involution $c_F$ has a unique fixed point on $\Sigma_P$, which 
corresponds to the intersection of the 
folding curve $F_P \cap \partial{D}$ of $c_F$ with the graph $\Sigma_P$. 
Maximal here
means that the number of edges of the graph of the map is maximal, i.e.
$3g$.

It would be interesting to compute the generating series  
$$
{\rm divide}(t):=\sum_{g\in \N}d(g)t^g=1+t+2t^2+8t^3+36t^4+ \cdots,
$$ 
where $d(g)$ denotes the 
number of simple, relative and  generic immersions the  interval $[0,1]$ 
with $g$ double points in the  disk $D$, counted up to homeomorphism in the
source and image. We have taken the coefficients up to the term $t^4$ of 
${\rm divide}(t)$  from the listings of
simple, relative, free or oriented divides by Masaharu Ishikawa [I].
We ask to compare the numbers $d(g)$ with the numbers $m^+(g)$ of maximal 
maps of genus $g$ with orientation reversing involution having $1$ fixed point
on its graph.

\noindent
{\bf 2. Involutions induced by 
$\pi$-rotations.}

Our next goal is to visualize the two involutions $c_F$ and $T \circ c_F$. 
We assume that the monodromy diffeomorphism $T$ is chosen 
in its isotopy class such that $T \circ c_F$
is an involution of the fiber surface $F_P$. In fact, we assume that
the monodromy $T$ was given by a monodromy vector field $X$, which satisfies
$X=X^c$.

\begin{figure}
\input{fig-5.pstex_t}
\center{Fig. $5$. {\rm Dehn twist as composition of two involutions.}}
\end{figure}

We assume for simplicity, that the divide $P$ meets the boundary of the disk
$D$.
The involution $c_F$ is an orientation reversing diffeomorphism of $F_1=F_P$,
which fixes pointwise the system of arcs 
$a:=F_P \cap \partial{D}$. The Lefschetz number
of the orientation reversing involution 
$T \circ c_F$ is equal to the Lefschetz number of $c_F$ by Th. $1$. 
Hence the 
involution $T \circ c_F$ also
fixes  pointwise exactly a system of  arcs $b$ on the fiber surface $F_P$ 
with $\chi(a)=\chi(b)$, i.e.
the systems $a$ and $b$ consist of the same number of arcs. 

Let 
$z:Z \to S^3 \setminus K_P$ be the infinite
cyclic covering  of the knot complement. Let $X^Z$ be the lift of the vector
field $X$. Let $T_{{1 \over 2}}:Z \to Z$ be the flow diffeomorphism of 
the vector field $X^Z$ with stopping time ${1 \over 2}$. 
Let $F_1' \subset Z$ be a lift of the fiber surface
$F_1$ and let $F_{-1}'$ be the lift $T_{{1 \over 2}}(F_1')$ 
of the fiber surface $F_{-1}$. The involution $c_F$ can first 
be lifted unambiguously to an 
involution of $F_1'$ and then extended unambiguously 
to an involution
$A$ of $Z$ such that $A$ maps the vector field $X^Z$ to its opposite. 
The involution $A$ is a $\pi$-rotation about a lift of the system of arcs 
$a$ into $F_1'$.
The involution $(x,u) \mapsto (x,-u)$ induces an involution $c_{F_{-1}}$ 
of the fiber surface $F_{-1}$ above $-1$, and,
as we have done for the involution $c_F$,  
the involution $c_{F_{-1}}$ lifts unambiguously to a 
$\pi$-rotation $B:Z \to Z$ about a system of  arcs  
$b'' \subset F_{-1}'$ that
satisfies
$z(b'')=b':=\partial{D} \cap F_{-1}$. 
We have $B \circ A=T_{{1 \over 2}}\circ T_{{1 \over 2}}=T_1$. Since 
$T_1$ is a lift of the monodromy, it follows from $(T \circ c_F) \circ c_F=T$, 
that
the involution $T \circ c_F$ fixes the system of arcs 
$b:=z(T_{{1 \over 2}}^{-1}(b''))$.

\begin{center}
\input{fig-6.pstex_t}
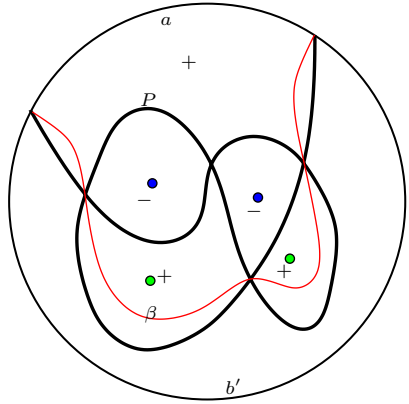
\center{Fig. $6$. Divide for the $E_8$ singularity and the curve $\beta$.}
\end{center}

In [AC4] we have computed the monodromy
diffeomorphism $T_+:F_1 \to F_{-1}$ for which $T_{{1 \over 2}}$ is a lift to
$Z$ as a product of half Dehn twists for a divide. 
From the above we deduce, that the
involution $T \circ c_F:F_1 \to F_1$ is the composition 
$T_+^{-1} \circ c_{F_{-1}} \circ T_+$. We see that $T \circ c_F$ fixes the 
arc $b=T_+^{-1}(\partial{D} \cap F_{-1})$. We also see that the involution
$c_F \circ T$ of $F_1$ fixes the system of arcs 
$c_F(b)=T_+(\partial{D} \cap F_{-1})$.
Both arc systems $b$ and $c_F(b)$ have equal  projections $\beta$ into $D$.
See Fig. $5,6$ for examples. 

The composition of the orientation reversing
involutions, that fix pointwise the arc systems $a,b$ of Fig. $5$  on the cylinder surface,
is indeed a Dehn twist. To see this, think of  the cylinder as 
$[0,\pi] \times S^1$; the two involutions induce on the circle 
$\{\theta\} \times S^1$ reflections about diameters that make an 
angle $\theta$, so the
composition is a rotation of angle $2\theta$ of the circle 
$\{\theta\} \times S^1$.

In Fig. $6$ the curve $\beta$ is drawn for a more complicated divide for the
plane curve singularity $E_8$ with equation $y^3-x^5=0$.

\noindent
{\bf References.}

\noindent
[AC1]
Norbert A'Campo,
{\it Le Groupe de Monodromie du D\'eploiement des Singularit\'es
Isol\'ees de Courbes Planes I},
Math. Ann. {\bf 213} (1975),
1--32.

\noindent
[AC2]
Norbert A'Campo,
{\it Le Groupe de Monodromie du D\'eploiement des Singularit\'es
Isol\'ees de Courbes Planes II}, 
Actes du Congr\`es International des Math\'e\-ma\-ti\-ciens, Vancouver (1974)
395--404.

\noindent
[AC3]
Norbert A'Campo,
{\it Real deformations and complex topology of plane curve singularities},
Annales de la Facult\'e des Sciences de Toulouse {\bf 8} (1999) 1, 5--23, 
erratum {\it ibid.} 2, 343.

\noindent
[AC4]
Norbert A'Campo,
{\it Generic immersions of curves, knots,
monodromy and gordian number},
Publ. Math. I.H.E.S. {\bf 88} (1998), 151-169, (1999).

\noindent
[AC5]
Norbert A'Campo,
{\it Planar trees, slalom curves and hyperbolic knots},
Publ. Math. I.H.E.S. {\bf 88} (1998), 171-180, (1999).

\noindent
[AC6]
Norbert A'Campo,
{\it Quadratic vanishing cycles, reduction curves and 
reduction of the monodromy
group of plane curve singularities}, 
Tohoku Math. Journal {\bf 53} (2001), 533--552.

\noindent
[B-K]
Ludwig Balke and Rainer Kaenders,
{\it On certain type of Coxeter-Dynkin diagrams of plane curve
singularities}, Topology {\bf 35} (1995), 39--54.

\noindent
[C]
J. Callahan,
{\it The double cusp has five minima},
Math. Proc. Cambridge Philos. Soc. 
{\bf 84} 3, (1978) 537--538. 

\noindent
[K]
Akio Kawauchi,
{\it A survey of knot theory},
Birkh\"auser, 
Basel 1996.

\noindent
[G-Z1]
S. M. Gusein-Zade,
{\it Matrices d'intersections pour certaines singularit\'es de
fonctions de 2 variables},
Funkcional. Anal. i Prilozen
{\bf 8} (1974),
11--15.

\noindent
[G-Z2]
S. M. Gusein-Zade,
{\it Index of a singular point of a gradient vector field},
Funkcional. Anal. i Prilozen
{\bf 18} (1984),
6--10.

\noindent
[H-To]
Wolfgang Heil and Jeffrey L.  Tollefson,
{\it Deforming homotopy involutions of $3$-manifolds to involutions},
Topology {\bf 17} 4, (1978) 349--365.

\noindent
[H-C]
D. Hilbert, S. Cohn-Vossen, {\it Anschauliche Geometrie}, Springer, Berlin, 1932.
Reprint: Dover, New York, 1944. English transl.: {\it Geometry and the Imagination}, Chelsea, New York, 1952.

\noindent
[I]
Masaharu Ishikawa,
{\it Links of plane curves, fibredness, quasipositivity and Jones polynomials},
Inauguraldissertation, Universit\"at Basel, 2001.

\noindent
[M]
J. Milnor,
{\it Singular points on complex hypersurfaces},
Ann. of Math. Studies, 
Princeton University Press, 1968.

\noindent
[S]
Mathias Schulze,
Email communications,  November 2001 and Thesis at Kaiserslautern in progress.

\noindent
[T]
B. Teissier,
{\it Un exemple de classe d'\'equisingularit\'e irrationelle}, 
C. R. Acad. Sci. Paris, {\bf 311} 2 (1990), 111--113.

\noindent
[To]
Jeffrey L. Tollefson, 
{\it Periodic homeomorphisms of $3$-manifolds fibered over $S\sp{1}$},
Trans. Amer. Math. Soc.,{\bf 223} (1976)
223--234, with Erratum Trans. Amer. Math. Soc., {\bf 243} (1978) 309--310.

\noindent
Universit\"at Basel, Rheinsprung 21, CH-4051  Basel.
\end{document}

%% file: fig-1.pstex_t
\begin{picture}(0,0)%
\includegraphics{fig-1.pstex}%
\end{picture}%
\setlength{\unitlength}{1243sp}%
\begingroup\makeatletter\ifx\SetFigFont\undefined%
\gdef\SetFigFont#1#2#3#4#5{%
  \reset@font\fontsize{#1}{#2pt}%
  \fontfamily{#3}\fontseries{#4}\fontshape{#5}%
  \selectfont}%
\fi\endgroup%
\begin{picture}(9590,9588)(1731,-8680)
\end{picture}

%% file: fig-2.pstex_t
\begin{picture}(0,0)%
\includegraphics{fig-2.pstex}%
\end{picture}%
\setlength{\unitlength}{1657sp}%
\begingroup\makeatletter\ifx\SetFigFont\undefined%
\gdef\SetFigFont#1#2#3#4#5{%
  \reset@font\fontsize{#1}{#2pt}%
  \fontfamily{#3}\fontseries{#4}\fontshape{#5}%
  \selectfont}%
\fi\endgroup%
\begin{picture}(8106,8106)(2248,-8614)
\put(6076,-2311){\makebox(0,0)[lb]{\smash{\SetFigFont{10}{12.0}{\rmdefault}{\mddefault}{\updefault}
$P$
}}}
\put(4906,-4696){\makebox(0,0)[lb]{\smash{\SetFigFont{10}{12.0}{\rmdefault}{\mddefault}{\updefault}
$S_P^+$
}}}
\put(5536,-7351){\makebox(0,0)[lb]{\smash{\SetFigFont{10}{12.0}{\rmdefault}{\mddefault}{\updefault}
$S_P^-$
}}}
\end{picture}

%% file: fig-3.pstex_t
\begin{picture}(0,0)%
\includegraphics{fig-3.pstex}%
\end{picture}%
\setlength{\unitlength}{1657sp}%
\begingroup\makeatletter\ifx\SetFigFont\undefined%
\gdef\SetFigFont#1#2#3#4#5{%
  \reset@font\fontsize{#1}{#2pt}%
  \fontfamily{#3}\fontseries{#4}\fontshape{#5}%
  \selectfont}%
\fi\endgroup%
\begin{picture}(12464,8593)(361,-8115)
\put(5671,-2851){\makebox(0,0)[lb]{\smash{\SetFigFont{10}{12.0}{\familydefault}{\mddefault}{\updefault}
$e_z$
}}}
\put(7201,-1636){\makebox(0,0)[lb]{\smash{\SetFigFont{10}{12.0}{\familydefault}{\mddefault}{\updefault}
$\delta_z$
}}}
\put(3376,-5056){\makebox(0,0)[lb]{\smash{\SetFigFont{10}{12.0}{\rmdefault}{\mddefault}{\updefault}
$+$
}}}
\put(5761,-6406){\makebox(0,0)[lb]{\smash{\SetFigFont{10}{12.0}{\rmdefault}{\mddefault}{\updefault}
$-$
}}}
\put(7156,-151){\makebox(0,0)[lb]{\smash{\SetFigFont{10}{12.0}{\rmdefault}{\mddefault}{\updefault}
$-$
}}}
\put(10486,-6046){\makebox(0,0)[lb]{\smash{\SetFigFont{10}{12.0}{\rmdefault}{\mddefault}{\updefault}
$+$
}}}
\put(3646,-6316){\makebox(0,0)[lb]{\smash{\SetFigFont{10}{12.0}{\rmdefault}{\mddefault}{\updefault}
$P$
}}}
\put(7966,-6676){\makebox(0,0)[lb]{\smash{\SetFigFont{10}{12.0}{\rmdefault}{\mddefault}{\updefault}
$P$
}}}
\put(6751,-3301){\makebox(0,0)[lb]{\smash{\SetFigFont{10}{12.0}{\rmdefault}{\mddefault}{\updefault}
$z$
}}}
\put(361,-2716){\makebox(0,0)[lb]{\smash{\SetFigFont{10}{12.0}{\rmdefault}{\mddefault}{\updefault}
$maximum$
}}}
\end{picture}

%% file: fig-4.pstex_t
\begin{picture}(0,0)%
\includegraphics{fig-4.pstex}%
\end{picture}%
\setlength{\unitlength}{1657sp}%
\begingroup\makeatletter\ifx\SetFigFont\undefined%
\gdef\SetFigFont#1#2#3#4#5{%
  \reset@font\fontsize{#1}{#2pt}%
  \fontfamily{#3}\fontseries{#4}\fontshape{#5}%
  \selectfont}%
\fi\endgroup%
\begin{picture}(11788,8751)(92,-8273)
\put(3691,-6361){\makebox(0,0)[lb]{\smash{\SetFigFont{10}{12.0}{\rmdefault}{\mddefault}{\updefault}
+
}}}
\put(1801,-3931){\makebox(0,0)[lb]{\smash{\SetFigFont{10}{12.0}{\rmdefault}{\mddefault}{\updefault}
+
}}}
\put(4456, 29){\makebox(0,0)[lb]{\smash{\SetFigFont{10}{12.0}{\rmdefault}{\mddefault}{\updefault}
$P$
}}}
\put(1306,-7396){\makebox(0,0)[lb]{\smash{\SetFigFont{10}{12.0}{\rmdefault}{\mddefault}{\updefault}
$P$
}}}
\put(10621,-5281){\makebox(0,0)[lb]{\smash{\SetFigFont{10}{12.0}{\rmdefault}{\mddefault}{\updefault}
$P$
}}}
\put(1936,-2581){\makebox(0,0)[lb]{\smash{\SetFigFont{10}{12.0}{\rmdefault}{\mddefault}{\updefault}
$maximum$
}}}
\put(4726,-4426){\makebox(0,0)[lb]{\smash{\SetFigFont{10}{12.0}{\rmdefault}{\mddefault}{\updefault}
$minimum$
}}}
\put(5311,-3571){\makebox(0,0)[lb]{\smash{\SetFigFont{10}{12.0}{\rmdefault}{\mddefault}{\updefault}
$-$
}}}
\put(6976,-3481){\makebox(0,0)[lb]{\smash{\SetFigFont{10}{12.0}{\rmdefault}{\mddefault}{\updefault}
$+$
}}}
\put(3511,-3571){\makebox(0,0)[lb]{\smash{\SetFigFont{10}{12.0}{\rmdefault}{\mddefault}{\updefault}
$+$
}}}
\put(4951,-5956){\makebox(0,0)[lb]{\smash{\SetFigFont{10}{12.0}{\rmdefault}{\mddefault}{\updefault}
$+$
}}}
\put(9451,-6001){\makebox(0,0)[lb]{\smash{\SetFigFont{10}{12.0}{\rmdefault}{\mddefault}{\updefault}
$-$
}}}
\put(5266,-241){\makebox(0,0)[lb]{\smash{\SetFigFont{10}{12.0}{\rmdefault}{\mddefault}{\updefault}
$-$
}}}
\put(901,-6091){\makebox(0,0)[lb]{\smash{\SetFigFont{10}{12.0}{\rmdefault}{\mddefault}{\updefault}
$-$
}}}
\put(4411,-7216){\makebox(0,0)[lb]{\smash{\SetFigFont{10}{12.0}{\rmdefault}{\mddefault}{\updefault}
$maximum$
}}}
\put(6841,-2716){\makebox(0,0)[lb]{\smash{\SetFigFont{10}{12.0}{\rmdefault}{\mddefault}{\updefault}
$maximum$
}}}
\end{picture}

%% file: fig-5.pstex_t
\begin{picture}(0,0)%
\includegraphics{fig-5.pstex}%
\end{picture}%
\setlength{\unitlength}{1450sp}%
\begingroup\makeatletter\ifx\SetFigFont\undefined%
\gdef\SetFigFont#1#2#3#4#5{%
  \reset@font\fontsize{#1}{#2pt}%
  \fontfamily{#3}\fontseries{#4}\fontshape{#5}%
  \selectfont}%
\fi\endgroup%
\begin{picture}(10518,5621)(1318,-6559)
\put(4141,-2581){\makebox(0,0)[lb]{\smash{\SetFigFont{9}{10.8}{\rmdefault}{\mddefault}{\updefault}\special{ps: gsave 0 0 0 setrgbcolor}$+$\special{ps: grestore}}}}
\put(2251,-4111){\makebox(0,0)[lb]{\smash{\SetFigFont{9}{10.8}{\rmdefault}{\mddefault}{\updefault}\special{ps: gsave 0 0 0 setrgbcolor}$+$\special{ps: grestore}}}}
\put(5086,-1951){\makebox(0,0)[lb]{\smash{\SetFigFont{9}{10.8}{\rmdefault}{\mddefault}{\updefault}\special{ps: gsave 0 0 0 setrgbcolor}$a$\special{ps: grestore}}}}
\put(1486,-4966){\makebox(0,0)[lb]{\smash{\SetFigFont{9}{10.8}{\rmdefault}{\mddefault}{\updefault}\special{ps: gsave 0 0 0 setrgbcolor}$a$\special{ps: grestore}}}}
\put(1801,-3211){\makebox(0,0)[lb]{\smash{\SetFigFont{9}{10.8}{\rmdefault}{\mddefault}{\updefault}\special{ps: gsave 0 0 0 setrgbcolor}$P$\special{ps: grestore}}}}
\put(6571,-3346){\makebox(0,0)[lb]{\smash{\SetFigFont{9}{10.8}{\rmdefault}{\mddefault}{\updefault}\special{ps: gsave 0 0 0 setrgbcolor}$a$\special{ps: grestore}}}}
\put(11836,-3346){\makebox(0,0)[lb]{\smash{\SetFigFont{9}{10.8}{\rmdefault}{\mddefault}{\updefault}\special{ps: gsave 0 0 0 setrgbcolor}$a$\special{ps: grestore}}}}
\put(9946,-3751){\makebox(0,0)[lb]{\smash{\SetFigFont{9}{10.8}{\rmdefault}{\mddefault}{\updefault}\special{ps: gsave 0 0 0 setrgbcolor}$b$\special{ps: grestore}}}}
\put(9901,-2941){\makebox(0,0)[lb]{\smash{\SetFigFont{9}{10.8}{\rmdefault}{\mddefault}{\updefault}\special{ps: gsave 0 0 0 setrgbcolor}$b$\special{ps: grestore}}}}
\put(3601,-2086){\makebox(0,0)[lb]{\smash{\SetFigFont{9}{10.8}{\rmdefault}{\mddefault}{\updefault}\special{ps: gsave 0 0 0 setrgbcolor}$\beta$\special{ps: grestore}}}}
\put(4276,-3121){\makebox(0,0)[lb]{\smash{\SetFigFont{9}{10.8}{\rmdefault}{\mddefault}{\updefault}\special{ps: gsave 0 0 0 setrgbcolor}$\beta$\special{ps: grestore}}}}
\put(4456,-4786){\makebox(0,0)[lb]{\smash{\SetFigFont{9}{10.8}{\rmdefault}{\mddefault}{\updefault}\special{ps: gsave 0 0 0 setrgbcolor}$b'$\special{ps: grestore}}}}
\put(1891,-2086){\makebox(0,0)[lb]{\smash{\SetFigFont{9}{10.8}{\rmdefault}{\mddefault}{\updefault}\special{ps: gsave 0 0 0 setrgbcolor}$b'$\special{ps: grestore}}}}
\put(1621,-6451){\makebox(0,0)[lb]{\smash{\SetFigFont{9}{10.8}{\rmdefault}{\mddefault}{\updefault}\special{ps: gsave 0 0 0 setrgbcolor}${\rm Divide}$ for $xy=0$\special{ps: grestore}}}}
\put(6571,-6451){\makebox(0,0)[lb]{\smash{\SetFigFont{9}{10.8}{\rmdefault}{\mddefault}{\updefault}\special{ps: gsave 0 0 0 setrgbcolor}Fiber surface with arc systems $ a,b$\special{ps: grestore}}}}
\end{picture}

%% file: fig-6.pstex_t
\begin{picture}(0,0)%
\includegraphics{fig-6.pstex}%
\end{picture}%
\setlength{\unitlength}{1119sp}%
\begingroup\makeatletter\ifx\SetFigFont\undefined%
\gdef\SetFigFont#1#2#3#4#5{%
  \reset@font\fontsize{#1}{#2pt}%
  \fontfamily{#3}\fontseries{#4}\fontshape{#5}%
  \selectfont}%
\fi\endgroup%
\begin{picture}(8836,8834)(1883,-8348)
\put(5266, 29){\makebox(0,0)[lb]{\smash{\SetFigFont{7}{8.4}{\rmdefault}{\mddefault}{\updefault}
$a$
}}}
\put(6706,-8161){\makebox(0,0)[lb]{\smash{\SetFigFont{7}{8.4}{\rmdefault}{\mddefault}{\updefault}
$b'$
}}}
\put(5716,-916){\makebox(0,0)[lb]{\smash{\SetFigFont{7}{8.4}{\rmdefault}{\mddefault}{\updefault}
$+$
}}}
\put(4816,-1771){\makebox(0,0)[lb]{\smash{\SetFigFont{7}{8.4}{\rmdefault}{\mddefault}{\updefault}
$P$
}}}
\put(4906,-6496){\makebox(0,0)[lb]{\smash{\SetFigFont{7}{8.4}{\rmdefault}{\mddefault}{\updefault}
$\beta$
}}}
\put(5176,-5641){\makebox(0,0)[lb]{\smash{\SetFigFont{7}{8.4}{\rmdefault}{\mddefault}{\updefault}
$+$
}}}
\put(4726,-3976){\makebox(0,0)[lb]{\smash{\SetFigFont{7}{8.4}{\rmdefault}{\mddefault}{\updefault}
$-$
}}}
\put(7156,-4201){\makebox(0,0)[lb]{\smash{\SetFigFont{7}{8.4}{\rmdefault}{\mddefault}{\updefault}
$-$
}}}
\put(7831,-5551){\makebox(0,0)[lb]{\smash{\SetFigFont{7}{8.4}{\rmdefault}{\mddefault}{\updefault}
$+$
}}}
\end{picture}